\documentclass[12pt]{amsart}
\usepackage[colorlinks=true,citecolor=black,linkcolor=black,urlcolor=blue]{hyperref}
\usepackage{amsmath,pifont,enumitem,mathtools}
\usepackage[english,  activeacute]{babel}
\usepackage[utf8]{inputenc}
\usepackage[dvipsnames]{xcolor}
\usepackage{amssymb}
\usepackage{amsthm}
\usepackage{graphics,graphicx}
\usepackage{enumerate}
\usepackage{array}
\usepackage{bm}
\usepackage{a4wide}
\usepackage{color, url}
\usepackage{float}
\usepackage{xcolor}

\newcommand{\seqnum}[1]{\href{http://oeis.org/#1}{\underline{#1}}}

\usepackage[colorinlistoftodos]{todonotes}

\theoremstyle{plain}
\newtheorem{theorem}{Theorem}[section]

\newtheorem{coro}[theorem]{Corollary}

\theoremstyle{definition}

\newtheorem{openproblem}{Open Problem}[section]

\newcommand{\Cat}{{\mathcal{C}}}

\newcommand\scalemath[2]{\scalebox{#1}{\mbox{\ensuremath{\displaystyle #2}}}}

\def\val{{\textsf{val}}}
\def\wruns{{\textsf{wruns}}}
\def\wdruns{{\overline{\textsf{wruns}}}}
\def\sval{{\textsf{sval}}}
\def\lval{{\textsf{$\ell$-val}}}
\def\lpea{{\textsf{$\ell$-peak}}}
\def\peak{{\textsf{peak}}}
\def\runs{{\textsf{runs}}}
\def\druns{{\overline{\textsf{runs}}}}
\def\vsym{{\textsf{symv}}}
\def\psym{{\textsf{symp}}}

\newcommand{\tx}[1]{\texttt{#1}}

\title{Enumerating runs, valleys, and peaks in Catalan Words}

\date{\today}

\begin{document}

\author[J.-L Baril]{Jean-Luc Baril}
\address[J.-L Baril]{LIB, Universit\'e de Bourgogne Franche-Comt\'e,   B.P. 47 870, 21078, Dijon Cedex, France}
\email{barjl@u-bourgogne.fr}

\author[P.~E.~Harris]{Pamela E. Harris}
\author[K.~J.~Harry]{Kimberly J. Harry}
\author[M. M{c}Clinton]{Matt McClinton}

\address[P.~E.~Harris, K.~J.~Harry, M. M{c}Clinton]{Department of Mathematical Sciences, University of Wisconsin-Milwaukee, Milwaukee, WI 53211 United States}
\email{peharris@uwm.edu}\email{kjharry@uwm.edu}\email{mcclin33@uwm.edu}

\author[J. L. Ram\'{\i}rez]{Jos\'e L. Ram\'{\i}rez}
\address[J. L. Ram\'{\i}rez]{Departamento de Matem\'aticas,  Universidad Nacional de Colombia,  Bogot\'a, Colombia}
\email{jlramirezr@unal.edu.co}

\newcommand{\matt}[1]{\mbox{}{\sf\color{violet}[Matt: #1]}\marginpar{\color{violet}\Large$*$}}

\newcommand{\pamela}[1]{\textcolor{red}{Pamela: #1}}
\newcommand{\tr}[1]{\textcolor{red}{#1}}

\newcommand{\kim}[1]{\mbox{}{\sf\color{olive}[Kim: #1]}\marginpar{\color{olive}\Large$*$}} 

\newcommand{\jluc}[1]{\mbox{}{\sf\color{magenta}[J-Luc: #1]}\marginpar{\color{magenta}\Large$*$}}

\newcommand{\jose}[1]{\mbox{}{\sf\color{blue}[Jose: #1]}\marginpar{\color{blue}\Large$*$}}

\begin{abstract}
We provide generating functions, formulas, and asymptotic  expressions for the number of Catalan words based on the number of runs of ascents (descents), runs of weak ascents (descents), $\ell$-valleys, valleys, symmetric valleys, $\ell$-peaks, peaks, and symmetric peaks. We also establish some bijections with restricted Dyck paths and ordered trees that transports some statistics.
\end{abstract}

\maketitle

\section{Introduction}

A word $w=w_1w_2\cdots w_n$ over the set of nonnegative integers is called a \emph{Catalan word} if $w_1=\tx{0}$ and $\tx{0}\leq w_i\leq w_{i-1}+\tx{1}$ for $i=2, \dots, n$.  We let $|w|$ denote the length of~$w$. We denote by $\epsilon$ the \emph{empty word}, that is the unique word of length zero. For $n\geq 0$, let $\Cat_n$ denote the set of Catalan words of length $n$. 
We set  $\Cat\coloneq\bigcup _{n\geq0}\Cat_n$ and $\Cat^+\coloneq\bigcup _{n\geq1}\Cat_n$ be the set of nonempty Catalan words. 
For example,
\begin{align*}
\Cat_4=\left\{\begin{matrix} \texttt{0000},  \texttt{0001},  \texttt{0010},  \texttt{0011},  \texttt{0012},  \texttt{0100},  \texttt{0101},  \\ \texttt{0110},  \texttt{0111},  \texttt{0112},  \texttt{0120},  \texttt{0121},  \texttt{0122},  \texttt{0123}\end{matrix}\right\}.
\end{align*}
As may be expected given their name, the cardinality of the set $\Cat_n$ is given by the $n$th Catalan number $c_n=\frac{1}{n+1}\binom{2n}{n}$, see  \cite[Exercise 80]{Stanley2}.  
Catalan words have  been  previously studied within the framework of exhaustively
generation of Gray codes for growth-restricted words \cite{ManVaj}. Baril  et al. \cite{BGR, Baril2, Baril}   have investigated the distribution of descents and the last symbol in restricted Catalan words that avoid  a classical pattern or a pair of patterns of length at most three. Analogous results have been found  in \cite{Baril3, Baril4, AlejaRam} for restricted Catalan words avoiding a consecutive  pattern  or an ordered pair of relations of length at most three.   On the other hand, Callan et al. \cite{CallManRam} initiate the enumeration of statistics, including area and perimeter, on the polyominoes associated with Catalan words.   Additionally, we refer to \cite{Baril5, ManRamF, ManRamM, Toc}, where the authors investigate several combinatorial statistics on the polyominoes linked to both  Catalan  and Motzkin words. Recently, Shattuck \cite{Shattuck} began studying the number of occurrences of  different subwords of length at most three contained within Catalan words, such as descents, ascents, levels, etc. Our paper continues the study of different statistics or parameters in Catalan words (runs of ascents, run of weak ascents, $\ell$-valleys, $\ell$ peaks,...) and their connections with other combinatorial objects, like trees and  restricted Dyck paths. Notice that the joint distribution of patterns $121$ and $231$  studied in \cite{Shattuck} corresponds to our study for $\ell$-peaks whenever $\ell=1$. Similarly, the joint distribution of consecutive patterns $212$ and $312$  studied in \cite{Shattuck} corresponds to our study for $\ell$-valleys for $\ell=1$.

Let $w=w_1w_2\cdots w_n\in \Cat_n$.
As usual, we say that $w$ has an \emph{ascent} (\emph{descent}) at position $\ell$ if $w_\ell < w_{\ell+1}$ ($w_\ell > w_{\ell+1}$), where $\ell \in [n-1]$. Similarly, we define \emph{weak ascent} (resp. \emph{weak descent}) whenever the inequality is not strict. A \emph{run} (resp. \emph{weak run}) of ascents (resp. \emph{weak ascents}) in a word $w$ is a maximal subword of consecutive ascents (resp. weak ascents). The number of runs in $w$ is denoted by $\runs(w)$, and the number of weak runs in $w$ is denoted by $\wruns(w)$. The runs of descents and weak descents are defined similarly, and the  statistics will be  denoted $\druns(w)$ and $\wdruns(w)$, respectively. 
An $\ell$-\emph{valley}  in a Catalan word $w$ is a subword of the form $ab^\ell(b+1)$, where $a>b$ and  $\ell$ is a positive integer and $b^\ell$ denotes $\ell$ consecutive copies of the letter $b$.  If $\ell=1$, we say that it is a  \emph{short valley}.  
The number of $\ell$-valleys of $w$ is denoted by $\lval(w)$ and the number of  all $\ell$-valleys for $\ell\geq 1$ of $w$ is denoted by $\val(w)$.  
A \emph{symmetric  valley} is a valley of the form $a(a-1)^\ell a$ with $\ell\geq 1$. The number of symmetric valleys of $w$ is denoted by $\vsym(w)$. Analogously, we define the peak statistic. An $\ell$-\emph{peak} in $w$ is a subword of the form $a(a+1)^\ell b$, where $a\geq b$ and  $\ell$ is a positive integer.
The number of $\ell$-peaks of $w$ is denoted by $\lpea(w)$ and
the sum of all $\ell$-peaks for $\ell\geq 1$ of $w$ is denoted by $\peak(w)$. 
If $\ell=1$, we say that it is a  \emph{short peak}; and if $a=b$,  it is called a  \emph{symmetric peak}. The number of symmetric peaks of $w$ is denoted by $\psym(w)$.

Among our contributions, we give generating functions, formulas, and asymptotic  expressions for the number of Catalan words based on the number of runs of ascents (descents), runs of weak ascents (descent), $\ell$-valleys, valleys, symmetric valleys, $\ell$-peaks, peaks, and symmetric peaks. We also establish some bijections with restricted Dyck paths and ordered trees.  These statistics have been studied in detail for Dyck paths, see for example \cite{Deut, Man, MansourD, Sapounakis, Sapounakis2, Sapounakis3}.

We aggregate our results and the notation used in Table~\ref{tab:notation}.
\begin{table}[H]
\begin{center}
\resizebox{\textwidth}{!}{
\begin{tabular}{|l|c|c|c|c|c|c|}\cline{2-7}
 \multicolumn{1}{l|}{} 	& \multicolumn{6}{c|}{Statistics} \\ \cline{2-7}
 \multicolumn{1}{l|}{} 		& runs  of asc. 		& runs  of w. asc. &runs  of desc. 	&  runs  of w. desc.&$\ell$-valleys	& short valleys   \\ \hline
 statistic on $w$ & $\runs(w)$ & $\wruns(w)$ &	$\druns(w)$&$\wdruns(w)$ & $\lval(w)$ & 1-$\val(w)$  \\ 
bivariate g. function & $R(x,y)$ & $W(x,y)$ &$\bar{R}(x,y)$& $\bar{W}(x,y)$&$V_\ell(x,y)$ & $V_1(x,y)$ \\
distribution & $r(n,k)$ & $w(n,k)$ &$\bar{r}(n,k)$&$\bar{w}(n,k)$& $v_\ell(n,k)$ & $v_1(n,k)$ \\
total number over $\Cat_n$ & $r(n)$  &  $w(n)$ &$\bar{r}(n)$ &$\bar{w}(n)$& $v_\ell(n)$ & $v_1(n)$  \\
\hline
\multicolumn{1}{l|}{} 	&valleys &sym. valleys & $\ell$-peaks & short peaks & peaks & sym.  peaks  \\
\hline
statistic on $w$&$\val(w)$& $\vsym(w)$& $\lpea(w)$ & 1-$\peak(w)$& $\peak(w)$ & $\psym(w)$  \\
bivariate g. function &$V(x,y)$& $S(x,y)$ & $P_\ell(x,y)$  & $P_1(x,y)$ & $P(x,y)$ & $T(x,y)$ \\
distribution &  $v(n,k)$&$s(n,k)$ & $p_\ell(n,k)$ &  $p_1(n,k)$ & $p(n,k)$&$t(n,k)$\\
total number over $\Cat_n$&  $v(n)$& $s(n)$  &  $p_\ell(n)$ & $p_1(n)$ & $p(n)$ & $t(n)$ \\
\hline
\end{tabular}
}
\end{center}
\caption{Summary of notation and results for statistics considered.} \label{tab:notation}
\end{table}

\section{Basic Definitions}\label{sec:Catalan words}

\begin{figure}[H]
\centering
\includegraphics[scale=0.8]{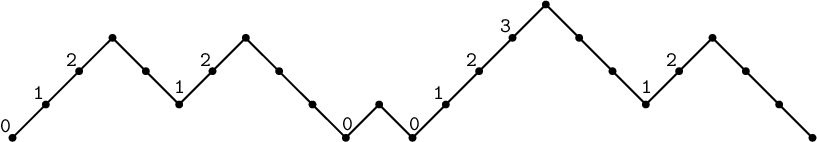}
\caption{Dyck path of the Catalan word \texttt{012120012312}.} \label{fig1}
\end{figure}

A {\it Dyck path} of semilength $n$ is a lattice path from $(0,0)$ to $(2n,0)$ consisting of $(1,1)$ up steps, and $(1,-1)$ down steps, in which the path never falls below the $x$-axis. Let $\mathcal{D}_n$ denote the set of all Dyck paths of semilength $n$.
Given a Dyck path $D\in \mathcal{D}_n$, we can associate a Catalan word in $\Cat_n$ formed by the $y$-coordinate of each initial point of the up steps in $D$. It is known that this construction yields a bijection between $\mathcal{D}_n$ and $\Cat_n$, see  \cite{Stanley2, CallManRam}.
For example, in Figure~\ref{fig1}, we illustrate the Dyck path associated to the Catalan word $\texttt{012120012312}\in \Cat_{12}$.

Throughout, we use the \emph{first return decomposition} of a Catalan word $w$, which~is \[w=\texttt{0}(w'+1)w'',\] where $w'$ and $w''$ are Catalan words ($w'$ and $w''$ could be empty), and where ($w'+1$) is the word obtained from $w'$ by adding $1$ at all these symbols. Note that whenever $w'$ is the empty word, denoted by $\epsilon$, then $(w'+1)$ remains the empty word. 
For example, the first return decomposition of $w=\texttt{012120012312}$ is given by setting $w'=\texttt{0101}$ and $w''=\texttt{0012312}$.  
For this word $w$, we have
$\runs(w)=5$, $\wruns(w)=4$, $\druns(w)=9$, $\wdruns(w)=8$, 1-$\val(w)=2$, 2-$\val(w)=1$, $\lval(w)=0$ $(\ell>2)$, $\vsym(w)=1$, 1-$\peak(w)=3$, $\lpea(w)=0$  $(\ell>1)$, and $\psym(w)=1$. 

\begin{figure}[H]
\centering
\includegraphics[scale=1]{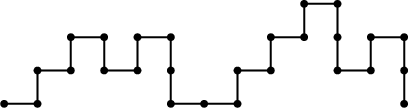}
\caption{Lattice diagram of the word \texttt{012120012312}.} \label{fig1b}
\end{figure}

Drawing Catalan words as lattice diagrams  on the plane proves to be a convenient representation. These diagrams are constructed using unit up steps $(0, 1)$, down steps $(0,-1)$, and horizontal steps $(1,0)$. Each symbol $w_i$ of a Catalan word is represented by the segment in between the points $(i-1, w_i)$ and $(i, w_i)$, and the vertical steps are inserted to obtain a connected diagram.   For example, in Figure~\ref{fig1b}, we illustrate the lattice diagram associated to the Catalan word of Figure~\ref{fig1}. 

\section{The Distribution of Runs}

\subsection{Runs of Ascents}\label{subsec:runs}
In order to count nonempty Catalan  words according to the length and the number runs of ascents, we introduce the following bivariate generating function
$$R(x,y)=\sum_{w \in\Cat^+}x^{|w|}y^{\runs(w)}=\sum_{n\geq 1}x^{|w|}\sum_{w\in\Cat_n}y^{\runs(w)},$$ where the coefficient of $x^ny^k$ is the number of Catalan words of length $n$ with $k$ runs of  ascents.

In Theorem~\ref{teo1}, we give an expression for this generating function.
\begin{theorem}\label{teo1}
The generating function for the number of nonempty Catalan words with respect to the length and the number of runs of ascents is
$$R(x,y)=\frac{1 - x (1 - y) -\sqrt{(1 - x - x y)^2 - 4 x^2 y}}{2x}.$$
\end{theorem}
\begin{proof}
Let $w$ be a nonempty Catalan word, and let $w=\texttt{0}(w'+1)w''$ be the first return decomposition of $w$, with $w', w''\in \Cat$.  
If $w'=w''=\epsilon$, then $w=\texttt{0}$, and its generating function is $xy$.  
If $w''=\epsilon$ and $w'\neq \epsilon$, then $w=\texttt{0}(w'+1)$, and its generating function is $xR(x,y)$.  
Similarly,  if $w'=\epsilon$ and $w''\neq \epsilon$, then $w=\texttt{0}w''$, and then its
generating function is $xyR(x,y)$ because we have made an extra run.
Finally, if $w'$ and $w''$ are not empty, then $w=\texttt{0}(w'+1)w''$, and then its
generating function is given by $xR^2(x,y)$. Therefore, we have the functional equation
$$R(x,y)=xy + x(1+y)R(x,y)+xR^2(x,y).$$
Solving this equation, we obtain the desired result. 
\end{proof}

\begin{figure}[H]
\centering
\includegraphics[scale=0.8]{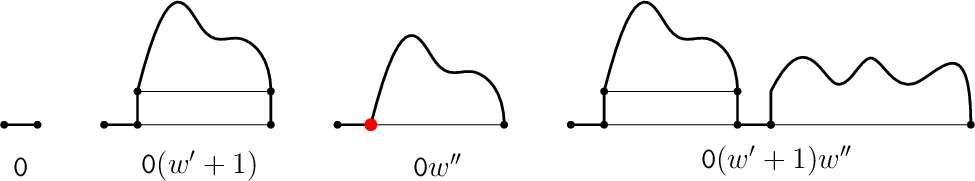}
\caption{Decomposition of a nonempty Catalan word.} \label{deco1}
\end{figure}

Figure~\ref{deco1}  illustrates the 
cases given in the decomposition described in the proof of Theorem~\ref{teo1}. The red point denotes the new run of ascents. 
 
Let $r(n,k)$ denote the number of Catalan words of length $n$ with exactly $k$  runs of  ascents, that is $r(n,k)=[x^ny^k]R(x,y)$, which denotes the coefficient of $x^ny^k$ in $R(x,y)$.
The first few values of this arrays are 
$$\mathcal{R}\coloneq[r(n,k)]_{n, k\geq 1}=
\begin{pmatrix}
 1 & 0 & 0 & 0 & 0 & 0 & 0 & 0 \\
 1 & 1 & 0 & 0 & 0 & 0 & 0 & 0 \\
 1 & 3 & 1 & 0 & 0 & 0 & 0 & 0 \\
 1 & \framebox{\textbf{6}} & 6 & 1 & 0 & 0 & 0 & 0 \\
 1 & 10 & 20 & 10 & 1 & 0 & 0 & 0 \\
 1 & 15 & 50 & 50 & 15 & 1 & 0 & 0 \\
 1 & 21 & 105 & 175 & 105 & 21 & 1 & 0 \\
 1 & 28 & 196 & 490 & 490 & 196 & 28 & 1 \\
\end{pmatrix}.$$

\begin{figure}[H]
\centering
\includegraphics[scale=0.8]{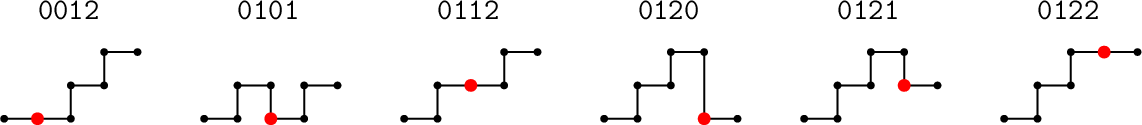}
\caption{Catalan words of length 4 with 2 runs of ascents.} \label{RunEx1}
\end{figure}

For example, $r(4,2)=6$, the entry boxed in $\mathcal{R}$ above, and the corresponding Catalan words (and lattice diagrams) are shown in Figure~\ref{RunEx1}.
The array $\mathcal{R}$ corresponds with the OEIS entry \cite[\seqnum{A001263}]{OEIS}, which is called the Narayana array.  Therefore, we have the combinatorial expression
$$r(n,k)=\frac{1}{k}\binom{n - 1}{k - 1} \binom{n}{k - 1}.$$
Notice, that this sequence counts the number of peaks in Dyck paths. 
Let $r(n)$ be the total number of runs of ascents over all Catalan words of length $n$.

\begin{coro}
For $n\geq 1$, we have
$r(n)=\binom{2n-1}{n}$.
\end{coro}
\begin{proof}
 By differentiating $R(x, y)$ with respect to $y$ and evaluating at $y= 1$, Theorem~\ref{teo1} gives
    \begin{align*}
 \left.\frac{\partial R(x,y)}{\partial y}\right|_{y=1}&=\frac{1+\sqrt{1 - 4x}}{2\sqrt{1 - 4x}}=\sum_{n\geq 1}\binom{2n-1}{n}x^n.
\end{align*}
By comparing the $n$th coefficient, we obtain the desired result. We use the generating  function of the central binomial coefficients:
\[\frac{1}{\sqrt{1-4x}}=\sum_{n\geq 0}\binom{2n}{n}x^n.\qedhere\]
\end{proof}

The first few values of the sequence $r(n)$ ($n\geq 1$) correspond with the OEIS entry \cite[\seqnum{A001700}]{OEIS}:
$$1, \quad 3, \quad 10, \quad 35, \quad 126, \quad 462, \quad 1716, \quad 6435, \quad 24310, \quad 92378, \dots.$$

Notice that in the notation of patterns in subwords, see  \cite{Shattuck}, the sequence $r(n)$ counts the total number of maximal patterns of the form $12\cdots (\ell-1)\ell$, with $\ell\geq 1$, in all Catalan words of length~$n$. 

\subsection{Runs of Weak Ascents} In order to count nonempty Catalan  words according to the length and the number runs of weak ascents, we introduce the following bivariate generating function $$W(x,y)=\sum_{w \in\Cat^+}x^{|w|}y^{\wruns(w)}=\sum_{n\geq 1}x^{|w|}\sum_{w\in\Cat_n}y^{\wruns(w)},$$ where the coefficient of $x^ny^k$ is the number of Catalan words of length $n$ with $k$ runs of weak ascents. 

In Theorem~\ref{teo2}, we give an expression for this generating function.

\begin{theorem}\label{teo2}
The generating function for the number of nonempty Catalan words with respect to the length and the number of runs of weak ascents is
$$W(x,y)=\frac{1 - 2 x -\sqrt{1 - 4 x + 4 x^2 - 4 x^2 y}}{2x}.$$
\end{theorem}

\begin{proof}
From a similar argument as in the proof of Theorem~\ref{teo1}, we obtain the functional equation
$$W(x,y)=xy + 2xW(x,y)+xW^2(x,y).$$
Solving this equation, we obtain the desired result.
\end{proof}

Let $w(n,k)$ denote the number of Catalan words of length $n$ with exactly $k$ runs of weak ascents, that is $w(n,k)=[x^ny^k]W(x,y)${, which denotes the coefficient of $x^ny^k$ in $W(x,y)$.}
The first few values of this array are 
$$\mathcal{W}\coloneq[w(n,k)]_{n, k\geq 1}=
\begin{pmatrix}
 1 & 0 & 0 & 0 & 0 \\
 2 & 0 & 0 & 0 & 0 \\
 4 & 1 & 0 & 0 & 0 \\
 8 & \framebox{\textbf{6}} & 0 & 0 & 0 \\
 16 & 24 & 2 & 0 & 0 \\
 32 & 80 & 20 & 0 & 0 \\
 64 & 240 & 120 & 5 & 0 \\
 128 & 672 & 560 & 70 & 0 \\
 256 & 1792 & 2240 & 560 & 14
\end{pmatrix}.$$

\begin{figure}[H]
\centering
\includegraphics[scale=0.8]{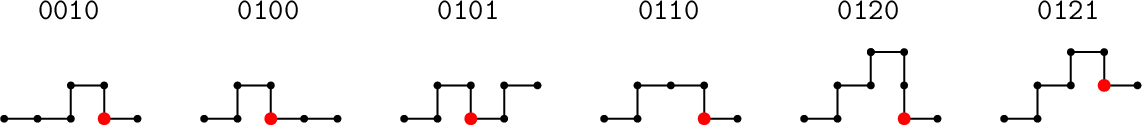}
\caption{Catalan words of length 4 with 2 runs of weak ascents.} \label{RunEx2}
\end{figure}

For example, $w(4,2)=6$, the entry boxed in $\mathcal{W}$ above, and the corresponding Catalan words (and lattice diagrams) are shown in Figure~\ref{RunEx2}.

\begin{theorem}\label{id1}
    For $n, k\geq 1$, we have
    $$w(n,k)=\frac{1}{n}\binom{n}{k}\binom{n-k}{k-1}2^{n - 2 k + 1}.$$
\end{theorem}
\begin{proof}
Rewriting the functional equation in the proof of Theorem~\ref{teo2}, we obtain  $$W(x,y) = xy + 2xW(x,y) + xW^2(x,y)= x\Phi(W(x,y)),$$ where $\Phi(z) = y+  2z + z^2$.
From the Lagrange inversion theorem (see for instance \cite{MSV}), we have
\begin{align*}
    w(n,k)&=[x^ny^k]W(x,y)=\frac{1}{n}[z^{n-1}y^k](y+2z+z^2)^n=\frac{1}{n}[z^{n-1}y^k]\sum_{i=0}^n\binom{n}{i}y^i(2z+z^2)^{n-i}\\
    &=\frac{1}{n}[z^{n-1}]\binom{n}{k}(2z+z^2)^{n-k}=\frac{1}{n}\binom{n}{k}[z^{n-1-(n-k)}]\sum_{i=0}^{n-k}\binom{n-k}{i}2^{n-k-i}z^i\\
    &=\frac{1}{n}\binom{n}{k}\binom{n-k}{k-1}2^{n - 2 k + 1}.\qedhere
\end{align*}
\end{proof}
From Theorem~\ref{id1}, we obtain the combinatorial identity for the Catalan numbers
\begin{align}\label{eq:cat identity}
c_n=\frac{1}{n}\sum_{k=1}^n \binom{n}{k}\binom{n-k}{k-1}2^{n-2k+1}.
\end{align}
\begin{openproblem}
Find a combinatorial interpretation for the identity in~\eqref{eq:cat identity}.
\end{openproblem}

The array $\mathcal{W}$ corresponds with the OEIS entry \cite[\seqnum{A091894}]{OEIS}. 
Among the combinatorial objects counted by this sequence, there are Dyck paths of semilength $n$ having exactly $k-1$ $UU$'s with midpoint at even height. 
Let $\mathcal{U}_n$ be the set of Dyck paths of semilength $n$ and $\mathcal{U}=\bigcup_{n\geq 0}\mathcal{U}_n$. For instance, there are $6$ Dyck paths of semilength $4$ having one occurrence $UU$ with midpoint at even height: 
\begin{align*}
U(UU)DUDDD, \quad U(UU)UDDDD,\quad  UDU(UU)DDD, \\
U(UU)DDDUD, \quad U(UU)DDUDD,\quad  UUD(UU)DDD.
\end{align*}

\begin{figure}[H]
\centering
\includegraphics[scale=0.8]{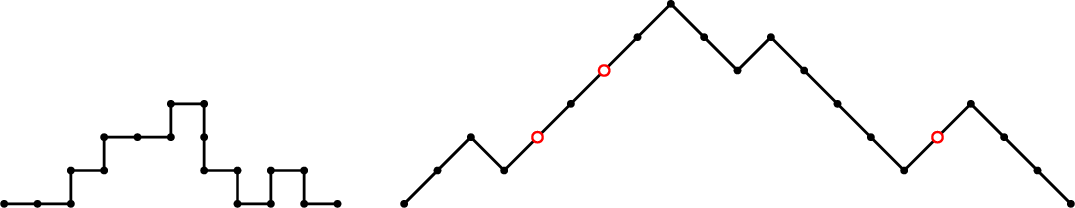}
\caption{An example for the bijection $\psi$ from $\mathcal{C}_{10}$ to  $\mathcal{U}_{10}$. } \label{bij1}
\end{figure}

We define recursively the bijection $\psi$ from $\mathcal{C}_n$ to  $\mathcal{U}_n$ as follows:\\
\scalemath{0.9}{$$\psi(w)=\begin{cases}
\epsilon & \mbox{if } w=\epsilon, \\
 UD\psi(u)& \mbox{if } w=\texttt{0}(u+1), \\
UU\psi(u_1) D U \cdots D U\psi(u_{a-1})DD \psi(u_a) & \mbox{if } w=\texttt{0}(u_1+1)\cdots \tx{0}(u_a+1),~  a\geq 2,~u_a\neq \epsilon,\\
UU\psi(u_1) D U \cdots D U\psi(u_{a-1})DU \psi(u_a)DD & \mbox{if } w=\tx{0}(u_1+1)\cdots \tx{0}(u_a+1)\tx{0},~ a\geq 1,\end{cases}$$}\\
where  $u_i$, $1\leq i\leq a$, are Catalan words (possibly empty except for the third case where $u_a\neq \epsilon $), and  where, as in the first return decomposition, $(u+1)$ corresponds to the word obtained from $u$ by increasing by one all letters of $u$. For instance, we have 
\begin{align*}
\psi(\tx{0012231010})&=\psi(\tx{0}\cdot \epsilon\cdot \tx{012231}\cdot\tx{01} \cdot \tx{0})\\
&=UU\psi(\epsilon)DU  \psi(\tx{01120})DU\psi(\tx{0}) DD\\
&= UU\cdot \epsilon\cdot DU\cdot UU \psi(\tx{001})DD \cdot DU \cdot UD \cdot DD\\
&=UUDUUU \cdot UU\psi(\epsilon)DD\psi(\tx{0})\cdot DDDUUDDD\\
&=UUDUUUUUDDUDDDDUUDDD.
\end{align*}
Refer to Figure~\ref{bij1}  for an illustration of the above example. The red points denote the midpoints at even height.

A simple observation allows us to prove that $\psi$ is a bijection. Now, let us check (by induction on the length) that $\psi$ transports the number of runs of weak ascents minus one into the number of occurrences $UU$ with midpoint at even height.  Let $\texttt{mid}(P)$ be the number of such occurrences in the path $P$.
By convenience, we set $\wruns(\epsilon)=1$, $\texttt{mid}(\epsilon)=0$ and the base case $n=0$ follows. Let $w$ satisfying $w=\texttt{0}(1+u_1)\cdots \tx{0}(1+u_a),~  a\geq 2,~u_a\neq \epsilon$, and let $k$ be  the number of nonempty $u_i$ for $1\leq i\leq a-1$. We have  
\begin{align*}\texttt{mid}(\psi(w))=&~\texttt{mid}(UU\psi(u_1) D U \cdots D U\psi(u_{a-1})DD \psi(u_a))\\
=&~k+\texttt{mid}(\psi(u_1))+ \cdots+\texttt{mid}(\psi(u_{a-1}))+\texttt{mid}(\psi(u_a)).
\end{align*}
Using the recurrence hypothesis, we have $\texttt{mid}(\psi(u_i))=\wruns(u_i)-1$ for $1\leq i\leq a$. So, we deduce
\begin{align*}
\texttt{mid}(\psi(w))=&~k+\big(\wruns(u_1)-1\big)+\cdots +\big(\wruns(u_{a-1})-1\big)+\big(\wruns(u_{a})-1\big)\\
=&~k+\wruns(u_1)+\cdots +\wruns(u_{a-1})-(a-1)+\wruns(u_{a})-1\\=&~ \wruns(\texttt{0}(1+u_1)\cdots \tx{0}(1+u_a))-1\\
=&~\wruns(w)-1.
\end{align*}
The second case where $w=\texttt{0}(1+u_1)\cdots \tx{0}(1+u_a)\tx{0},~ a\geq 1$ is handled {\it mutatis mutandis}.

Let $w(n)$ be the total number of runs of weak ascents over all Catalan words of length~$n$.

\begin{coro}
For $n\geq 1$, we have
$w(n)=\binom{2n-2}{n-1}$.
\end{coro}

The first few values of the sequence $w(n)$ corresponds to the OEIS entry \cite[\seqnum{A000984}]{OEIS}
$$1, \quad 2,  \quad 6, \quad 20, \quad 70, \quad 252, \quad 924, \quad 3432,\quad  12870, \dots.$$
Our combinatorial interpretation does not appear in the OEIS.

\subsection{Runs of Descents}
In order to count nonempty Catalan  words according to the length and the number runs of descents, we introduce the following bivariate generating function
$$\bar{R}(x,y)=\sum_{w \in\Cat^+}x^{|w|}y^{\druns(w)}=\sum_{n\geq 1}x^{|w|}\sum_{w\in\Cat_n}y^{\druns(w)},$$ where the coefficient of $x^ny^k$ is the number of Catalan words of length $n$ with $k$ runs of  descents.

In Theorem~\ref{teod1}, we give an expression for this generating function.
\begin{theorem}\label{teod1}
The generating function for the number of nonempty Catalan words with respect to the length and the number of runs of descents is
$$\bar{R}(x,y)=\frac{1-2 x y -\sqrt{4 x^{2} y^{2}-4 x^{2} y -4 x y +1}}{2 x}$$
\end{theorem}
\begin{proof}
Let $w$ be a nonempty Catalan word, and let $w=\texttt{0}(w'+1)w''$ be the first return decomposition of $w$, with $w', w''\in \Cat$.  
If $w'=w''=\epsilon$, then $w=\texttt{0}$, and its generating function is $xy$.  
If $w''=\epsilon$ and $w'\neq \epsilon$, then $w=\texttt{0}(w'+1)$, and then its generating function is $xy\bar{R}(x,y)$.   
Similarly,  if $w'=\epsilon$ and $w''\neq \epsilon$, then $w=\texttt{0}w''$, and then its generating function is $xy\bar{R}(x,y)$.
Finally, if $w'$ and $w''$ are not empty, then $w=\texttt{0}(w'+1)w''$, and then its generating function is given by $x\bar{R}^2(x,y)$. 
Therefore, we have the functional equation
$$\bar{R}(x,y)=xy + 2xy\bar{R}(x,y)+x\bar{R}^2(x,y).$$
Solving this equation yields the desired result. 
\end{proof}

Let $\bar{r}(n,k)$ denote the number of Catalan words of length $n$ with exactly $k$  runs of  descents, that is $\bar{r}(n,k)=[x^ny^k]\bar{R}(x,y)$, which denotes the coefficient of $x^ny^k$ in $\bar{R}(x,y)$.
The first few values of this arrays are 
$$\bar{\mathcal{R}}\coloneq[\bar{r}(n,k)]_{n, k\geq 1}=
\begin{pmatrix}
 1 & 0 & 0 & 0 & 0 & 0 & 0 & 0 \\
 0 & 2 & 0 & 0 & 0 & 0 & 0 & 0 \\
 0 & 0 & 6 & 8 & 0 & 0 & 0 & 0 \\
 0 & 0 & 2 & 24 & 16 & 0 & 0 & 0 \\
 0 & 0 & 0 & 20 & 80 & 32 & 0 & 0 \\
 0 & 0 & 0 & 5 & 120 & 240 & 64 & 0 \\
 0 & 0 & 0 & 0 & 70 & 560 & 672& 128 
\end{pmatrix}.$$
The array $\bar{\mathcal{R}}$ does not appear in the OEIS.
\begin{coro}\label{corod1}
For $n\geq 1$, we have
$\bar{r}(n)=\binom{2n}{n} - \binom{2n-2}{n-1}$.
\end{coro}
\begin{proof}
 By differentiating $\bar{R}(x, y)$ with respect to $y$ and evaluating at $y= 1$, Theorem~\ref{teod1} gives
    \begin{align*}
 \left.\frac{\partial \bar{R}(x,y)}{\partial y}\right|_{y=1}&=\frac{1-x-\sqrt{1-4 x }}{\sqrt{1-4 x }}.
\end{align*}
By comparing the $n$th coefficient, we obtain the desired result. We use the generating  function of the central binomial coefficients:
\[\frac{1}{\sqrt{1-4x}}=\sum_{n\geq 0}\binom{2n}{n}x^n.\qedhere\]
\end{proof}
The first few values of the sequence $\bar{r}(n)$ ($n\geq 1$) correspond with the OEIS entry \cite[\seqnum{A051924}]{OEIS}:
$$1, \quad 4, \quad 14, \quad 50, \quad 182, \quad 672, \quad 2508, \quad 9438, \quad 35750, \quad 136136, \dots .$$
\subsection{Runs of Weak Descents}
In a Catalan word of length $n$, the number of runs of ascents plus the number of runs of weak descents equals $n+1$. Due to the relation $r(n,k)=r(n,n+1-k)$ in the Narayana array (see Section~\ref{subsec:runs}), the number $\bar{w}(n,k)$ of Catalan words of length $n$ with $k$ runs of weak descents equals the number $r(n,k)$ of Catalan words of length $n$  with $k$ runs of ascents. Then, we refer to Section~\ref{subsec:runs}.

\section{The Distribution of Valleys}

\subsection{Valleys}

In order to count Catalan words according to the length and the number $\ell$-valleys, we introduce the following bivariate generating function $$V_\ell(x,y)=\sum_{w \in\Cat^+}x^{|w|}y^{\lval(w)}=\sum_{n\geq 1}x^{|w|}\sum_{w\in\Cat_n}y^{\lval(w)},$$
where $\lval(w)$ denotes the number of occurrences of subwords of the form $ab^\ell(b+1)$, and $a> b$, in $w$. 
The coefficient of $x^ny^k$ in $V_\ell (x,y)$ is the number of Catalan words of length $n$ with $k$ $\ell$-valleys.

In Theorem~\ref{teosv1l}, we give an expression for this generating function.

\begin{theorem}\label{teosv1l}
The generating function for nonempty  Catalan words with respect to the length and the  number of $\ell$-valleys is
$$V_\ell(x,y)=\frac{1 - 2 x - x^{\ell+1} (1 - y)-\sqrt{1 - 4 x + 2 x^{\ell+ 1}(1 - y) + x^{2 (\ell+1)} (1 - y)^2}}{2 x (1 - x^{\ell-1} (1 - y) + x^\ell (1 - y))}.
$$
\end{theorem}
\begin{proof} Let $w$ be a nonempty Catalan word, and let  $w=\texttt{0}(w'+1)w''$ be the first return decomposition, with $w', w''\in \Cat$. If $w'=w''=\epsilon$, then $w=\texttt{0}$, and its generating function is $x$.  If $w'\neq \epsilon$ and $w''=\epsilon$, then $w=\texttt{0}(w'+1)$, and its generating function is $xV_\ell (x,y)$.  
Similarly,  if $w'=\epsilon$ and $w''\neq \epsilon$, then $w=\texttt{0}w''$, and its generating function is $xV_\ell (x,y)$.  
Finally, if $w'\neq \epsilon$ and $w'' \neq \epsilon$, then $w=\texttt{0}(w'+1)w''$, and we distinguish two cases. 
If $w''$ is of the form $\tx{0}^{\ell-1}w'''$, where $w'''$ starts with $\tx{01}$, then $w=\texttt{0}(w'+1)\texttt{0}^{\ell-1}w'''$, and the generating function is 
$xyV_\ell(x,y)x^{\ell-1}A_\ell(x,y)$, 
where $A_\ell(x,y)=V_\ell(x,y)-x(1+V_\ell(x,y))$ is obtained using the complement of the generating function for the word $\texttt{0}$ and the words starting with $\texttt{00}$. Otherwise, the generating function is $xV_\ell (x,y)(V_\ell(x,y)-x^{\ell-1}A_\ell(x,y))$. 

Therefore, we have the functional equation
\begin{align*}V_\ell(x,y)=x + 2xV_\ell(x,y)&+xyV_\ell(x,y)x^{\ell-1}\left(V_\ell(x,y)-x(1+V_\ell(x,y)\right)+\\
&xV_\ell (x,y)\left(V_\ell(x,y)-x^{\ell-1}\left(V_\ell(x,y)-x(1+V_\ell(x,y))\right)\right).
\end{align*}
Solving this equation, we obtain the desired result.
\end{proof}

\begin{figure}[H]
\centering
\includegraphics[scale=0.8]{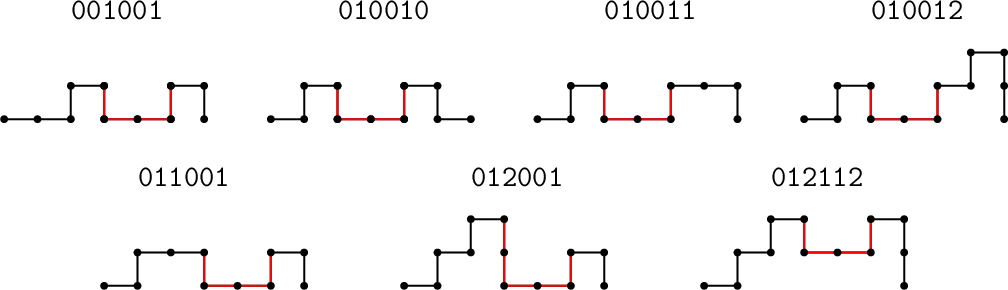}
\caption{Catalan words of length 6 with one $2$-valley.} \label{ExlVal}
\end{figure}

Let $v_\ell(n,k)$ denote the number of Catalan words of length $n$ with exactly $k$ $\ell$-valleys, that is $v_\ell(n,k)=[x^ny^k]V_\ell(x,y)$, which denotes the coefficient of $x^ny^k$ in $V_\ell(x,y)$.
For example, the first few values of this array for $\ell=2$ are 
$$\mathcal{V}_2\coloneq[v_2(n,k)]_{n\geq 1, k\geq 0}=
\begin{pmatrix}
 1 & 0 & 0 & 0 \\
 2 & 0 & 0 & 0 \\
 5 & 0 & 0 & 0 \\
 14 & 0 & 0 & 0 \\
 41 & 1 & 0 & 0 \\
 125 & \framebox{\textbf{7}} & 0 & 0 \\
 393 & 36 & 0 & 0 \\
 1266 & 163 & 1 & 0 \\
 4158 & 693 & 11 & 0 
\end{pmatrix}.$$
For example, $v_2(6,1)=7$, the entry boxed in $\mathcal{V}_2$ above, and the corresponding Catalan words (and lattice diagrams) are shown in Figure~\ref{ExlVal}. The red segments denote the valley.

Let $v_\ell(n)$ be the sum of all $\ell$-valleys in the set of  Catalan words of length $n$. 
\begin{coro}
The generating function of the sequence
$v_\ell(n)$  is 
$$\sum_{n\geq 1}v_\ell(n)x^n=\frac{x^{\ell-2}\left(1 - 5 x + 5 x^2 - (1 - 3 x + x^2)\sqrt{1-4x}\right)}{2\sqrt{1-4x}}.$$
Moreover, for all $n\geq 1$, we have
\begin{align}\label{eq:vell}
    v_\ell(n)=\binom{2(n - \ell) - 1}{n - \ell - 3},
\end{align}
 and an asymptotic approximation for $v_\ell(n)$ is $$\frac{2^{2(n-\ell)-1}}{\sqrt{\pi n}}.$$
\end{coro}
We compute asymptotics of the combinatorial sequences  using  classical singularity analysis (cf. \cite{flajolet}).

Taking $\ell=1$ in Theorem~\ref{teosv1l}, we obtain the generating function for nonempty Catalan words with respect to the length and the  number of short valleys
$$V_1(x,y)=\sum_{w \in\Cat^+}x^{|w|}y^{\text{1-}\val(w)}=\frac{1 - 2 x - x^2(1 - y)-\sqrt{1 - 4 x + 2 x^{2}(1 - y) + x^{4} (1 - y)^2}}{2 x (y + x(1 - y))}.
$$

Notice that in the notation of subwords patterns, see  \cite{Shattuck}, the sequence $v_1(n)$ counts the  joint distribution of the  patterns $212$ and $312$. 

Let $v_1(n,k)$ denote the number of Catalan words of length $n$ with exactly $k$ short valleys, that is $v_1(n,k)=[x^ny^k]V_1(x,y)$, which denotes the coefficient of $x^ny^k$ in $V_1(x,y)$. The first few values of this array are

$$\mathcal{V}_1=[v_1(n,k)]_{n\geq 1, k\geq 0}=
\begin{pmatrix}
 1 & 0 & 0 & 0 & 0\\
 2 & 0 & 0 & 0 & 0\\
 5 & 0 & 0 & 0 & 0\\
 13 & 1 & 0 & 0 & 0\\
 35 & \framebox{\textbf{7}} & 0 & 0& 0 \\
 97 & 34 & 1 & 0& 0 \\
 275 & 143 & 11 & 0 & 0\\
 794 & 558 & 77 & 1& 0 \\
 2327 & 2083 & 436 & 16 & 0
\end{pmatrix}.$$

\begin{figure}[H]
\centering
\includegraphics[scale=0.8]{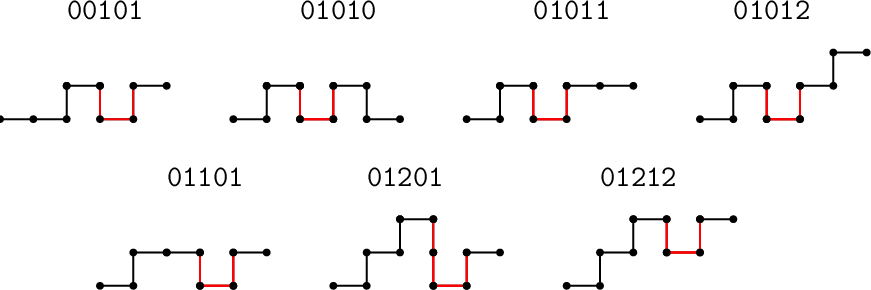}
\caption{Catalan words of length 5 with one short  valley.} \label{ExSV}
\end{figure}

For example, $v_1(5,1)=7$, the entry boxed in $\mathcal{V}_1$ above, and the corresponding Catalan words (and lattice diagrams) are shown in Figure~\ref{ExSV}.

The matrix $\mathcal{V}_1$ corresponds to the OEIS sequence \cite[\seqnum{A114492}]{OEIS}. Among the combinatorial objects counted by this sequence are Dyck paths of semilength $n$ having exactly $k$ occurrences of $DDUU$. For instance, the $7$  Dyck paths of semilength 5 having one short valley are  
\begin{center}
\begin{tabular}{lll}
$UDUU(DDUU)DD$, & $UUDU(DDUU)DD$,& $UU(DDUU)DDUD$,\\
  $UUUD(DDUU)DD$,&
  $UU(DDUU)UUUDDD$,&  $UU(DDUU)UUDUDD$,\\     $UUUD(DDUU)DDD$.&& \\
\end{tabular}
\end{center}
  
We recursively define  the bijection $\chi$ from $\mathcal{C}_n$ to  the set of Dyck paths of semilength $n$ as follows:
$$\chi(w)=\left\{\begin{array}{ll}
\epsilon & \mbox{if } w=\epsilon, \\
 U\chi(u)D& \mbox{if } w=\texttt{0}(1+u), \\
  UD\chi(u)& \mbox{if } w=\texttt{0}u, \\
U \chi(u) D\chi(v)& \mbox{if } w=\tx{0}(1+u)v,~u,v\neq \epsilon,
\end{array}\right.
$$
where $(1+u)$ corresponds to the word obtained from $u$ by increasing by one all letters of~$u$.
For instance, we have \begin{align*}
\chi(\tx{0012231010})&=\chi(\tx{0}\cdot \tx{012231010})\\
&=UD\chi(\tx{0}\cdot \tx{12231} \cdot \tx{010} )\\
&=
UDU\chi(\tx{0}\cdot\tx{112}\cdot \tx{0})D\chi(\tx{0}\cdot \tx{1}\cdot\tx{0}) \\
&=
UDU UUDUUDD D UDDUUDDUD.
\end{align*}

A simple observation allows us to prove that $\chi$ is a bijection. Now, let us check (by induction on the length) that $\chi$ transports the number of short valleys into the number of occurrences of $DDUU$.  Let $\texttt{dduu}(P)$ be the number of such occurrences in the path $P$. The  case $n=0$ is trivial. Let us consider the general decomposition of a Catalan word  $w=\tx{0}(1+u)v$. The cases where $u$ or $v$ are empty is easy to check, so we assume that $u$ and $v$ nonempty. We have  
$\texttt{dduu}(\chi(w))=~\texttt{dduu}(U\chi(u) D \chi(v))$.  We distinguish two cases.
 \begin{itemize}[leftmargin=.75in]
     \item[Case ($i$):] $v$ starts with $\tx{01}$. By definition, $\chi(v)$ starts with $UU$ and we have 
\begin{align*}
\texttt{dduu}(\psi(w))=& ~\texttt{dduu}(U\chi(u) D \chi(v))=~\texttt{dduu}(\chi(u))+1+\texttt{dduu}(\chi(v)).
\end{align*}
 Using the recurrence hypothesis, we deduce 
\begin{align*}\texttt{dduu}(\psi(w))=&~ 1+\sval(u)+\sval(v)=\sval(w).
\end{align*}
\item[Case $(ii)$:] $v$ does not start with $\tx{01}$. In the same way, we have 
 \begin{align*}
\texttt{dduu}(\psi(w))=&~ \texttt{dduu}(U\chi(u) D \chi(v))\\
=&~\texttt{dduu}(\chi(u))+\texttt{dduu}(\chi(v))=\sval(u)+\sval(v)=\sval(w).
\end{align*}
 \end{itemize}

Using a similar proof as for Theorem~\ref{teosv1l}, we generalize the result in order to obtain the following generating function for the
number Catalan words of length $n$ with respect to the number of valleys (we consider all $\ell$-valleys for $\ell\geq 1$).

\begin{theorem}\label{teov1}
The generating function for nonempty Catalan words with respect to the length and the  number of valleys is
$$V(x,y)=\frac{1 - 3 x + x^2 + x^2 y -\sqrt{(1 - 3 x + x^2 + x^2 y)^2-4 x^2 y (1 - x)^2}}{2xy(1 - x)}.$$
\end{theorem}

Let $v(n,k)$ denote the number of Catalan words of length $n$ with exactly $k$ valleys, that is $v(n,k)=[x^ny^k]V(x,y)$, which denotes the coefficient of $x^ny^k$ in $V(x,y)$.
The first few values of this arrays are 
$$\mathcal{V}=[v(n,k)]_{n\geq 1, k\geq 0}=
\begin{pmatrix}
 1 & 0 & 0 & 0 & 0 \\
 2 & 0 & 0 & 0 & 0 \\
 5 & 0 & 0 & 0 & 0 \\
 13 & 1 & 0 & 0 & 0 \\
 34 & \framebox{\textbf{8}} & 0 & 0 & 0 \\
 89 & 42 & 1 & 0 & 0 \\
 233 & 183 & 13 & 0 & 0 \\
 610 & 717 & 102 & 1 & 0 \\
 1597 & 2622 & 624 & 19 & 0 
\end{pmatrix}.$$

\begin{figure}[H]
\centering
\includegraphics[scale=0.8]{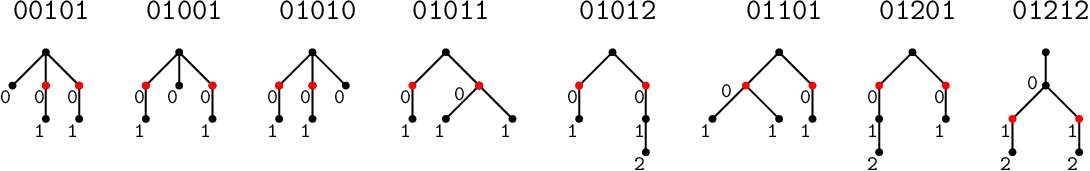}
\caption{Bijection $\psi$.} \label{BijTree}
\end{figure}

Note that $v(n,1)=F_{2n-1}$ ($n\geq 1$), where $F_n$ is the $n$th Fibonacci number.  
The array $\mathcal{V}$ corresponds to the OEIS entry \cite[\seqnum{A114502}]{OEIS}, where the $n$th term counts the number of ordered trees with $n$ edges and having exactly $k+1$ nodes all of whose children are leaves. 
For example, $v(5,1)=8$, the entry boxed in $\mathcal{V}$ above, and the corresponding Catalan words and trees are shown in Figure~\ref{BijTree}. 

\begin{theorem}
The set of ordered trees with $n$ edges and having exactly $k+1$ nodes all of whose children are leaves is in bijection with the set of Catalan words of length $n$ having $k$ valleys.
\end{theorem}
\begin{proof}
Recall that an ordered tree is either a single node, called root and denoted $()$ or equivalently $\bullet$, or a root whose some ordered subtrees $S_1, S_2, \ldots, S_k$, $k\geq 1$, denoted $(S_1,S_2,\ldots, S_k)$. A leave of an ordered tree is a node with no children. We define recursively the bijection $\psi$ from Catalan words of length $n$ to ordered trees with $n$ edges as follows: 
\begin{align*}
    (i)~ \psi(w)&=\bullet, \mbox{ if } w=\epsilon,\\
    (ii)~ \psi(w)&=(\psi(v_1),\psi(v_2),\ldots ,\psi(v_k)), \mbox{ if } w=\tx{0}(1+v_1)\cdots \tx{0}(1+v_k), k\geq 1, 1\leq i\leq k.
    \end{align*}
 Clearly, $\psi$ is a bijection, and the word $w$ corresponds to the level minus one  of the nodes (except the root) read in the preorder traversal of $\psi(w)$ (see Figure~\ref{BijTree}). Now let us prove by induction on the length that $\psi$ transports the number $\val(w)$ of valleys  plus one into the number $\texttt{nl}(\psi(w))$ of nodes whose all children are leaves in the tree $\psi(w)$. Assuming that $w=\tx{0}(1+v_1)\cdots \tx{0}(1+v_k)$, where $v_i$ are possibly empty,  we have
$$\texttt{nl}(\psi(w))=~\texttt{nl}\left((\psi(v_1),\psi(v_2),\ldots ,\psi(v_k))\right).$$

We distinguish two cases.
\begin{itemize}[leftmargin=.75in]
\item[Case ($i$):] For all $j$, $1\leq j\leq k$, we have  $v_j=\epsilon$. Thus, we have $w=\tx{0}^k$, $\psi(w)=(\bullet,\ldots, \bullet)$,  and  $\texttt{nl}(\psi(w))=1=\val(w)+1$.

    \item[Case ($ii$):] Otherwise, let $j_1,j_2, \ldots, j_a$, $a\geq 1$, the indices $j$ where $v_j\neq \epsilon$.
Thus, we have $\texttt{nl}(\psi(w))=\texttt{nl}(\psi(v_{j_1}))+\cdots +\texttt{nl}(\psi(v_{j_a}))$. Using the recurrence hypothesis, we obtain 
$\texttt{nl}(\psi(w))=\val(v_{j_1})+\cdots +\val(v_{j_a})+a=\val(w)+1$, which completes the induction.\qedhere
\end{itemize}
\end{proof}

Let $v(n)$ be the sum of all valleys in the set of  Catalan words of length $n$. 
\begin{coro}
The generating function of the sequence
$v(n)$  is 
$$\sum_{n\geq 0}v(n)x^n=\frac{1 - 5 x + 5 x^2 - (1 - 3 x + x^2)\sqrt{1-4x}}{2(1-x)x\sqrt{1-4x}}.$$ 
Moreover, for all $n\geq 1$, we have
\begin{align}\label{eq:sumvell}
    v(n)=\sum_{\ell=1}^{n-1}\binom{2(n - \ell) - 1}{n - \ell - 3},
\end{align}
and an asymptotic approximation for $v(n)$ is $$\frac{2^{2n}}{6\sqrt{\pi n}}.$$
\end{coro}

For $n\geq 4$, the first few values of the sequence $v(n)$ are 
$$1, \quad 8, \quad 44, \quad 209, \quad 924, \quad 3927, \quad  16303, \quad 66691,\quad  270181,\dots.$$
This sequence does not appear in the OEIS.

\begin{openproblem}
    Find a combinatorial proof for the equality in \eqref{eq:sumvell}. 
\end{openproblem}

\subsection{Symmetric Valleys}
A \emph{symmetric  valley} is a valley of the form $a(a-1)^\ell a$ with $\ell\geq 1$. Let $\vsym(w)$ denote  the number of  symmetric valleys in the word $w$. 
In order to count Catalan words according to the length and the number of symmetric valleys, we introduce the following bivariate generating function
generating function
$$S(x,y)=\sum_{w\in\Cat^+}x^{|w|}y^{\vsym(w)}=\sum_{n\geq 1}x^{|w|}\sum_{w\in\Cat_n}y^{\vsym(w)},$$
where the coefficient of $x^ny^k$ in $S(x,y)$ is the number of nonempty Catalan words of length~$n$ with $k$ symmetric $\ell$-valleys.

We remark that
the symmetric valley statistic has been recently studied in the context of $k$-ary words \cite{Asakly}, Dyck paths   \cite{FlorezJoseSymmDyck, Sergi},  and integer compositions \cite{Moreno}.

\begin{theorem}
    The generating function of the nonempty Catalan words with respect to the length and the number of symmetric valleys is 
    $$S(x,y)=\scalemath{0.85}{\frac{1 - 3 x + x^2 (3 - y) - 2 x^3 (1 - y) -\sqrt{1 - 6 x + x^2 (11 - 2 y) - 2 x^3 (5 - 3 y) + x^4 (5 - 6 y + y^2)}}{2 x (1 - x)  (1 - x (1 - y))}}.$$    
\end{theorem}
\begin{proof} Let $w$ be a nonempty Catalan word, and let  $w=\texttt{0}(w'+1)w''$ be the first return decomposition, with $w', w''\in \Cat$. If $w'=w''=\epsilon$, then $w=\texttt{0}$, and its generating function is $x$. If $w'\neq \epsilon$ and $w''=\epsilon$, then $w=\texttt{0}(w'+1)$, and its generating function is $xS(x,y)$. Similarly, if $w'=\epsilon$ and $w''\neq \epsilon$, then $w=\texttt{0}w''$ and its generating function is $xS(x,y)$.  Finally, if $w'$ and $w''$ are not empty, then $w=\texttt{0}(w'+1)w''$ and we have a new symmetric valley in the case where $w'$ ends with $\tx{0}$ and $w''\neq \tx{0}^\ell$ ($\ell\geq 1$).  The generating function for this case is
$$xS(x,y)^2 + x (xS(x,y) + x)y\left(S(x,y) - \frac{x}{1-x}\right) - x (xS(x,y) + x)\left(S(x,y) - \frac{x}{1-x}\right).$$
Summing  all the cases and solving the obtained functional equation yields the desired result.
\end{proof}

Let $s(n,k)$ denote the number of Catalan words of length $n$ with exactly $k$ symmetric valleys, that is $s(n,k)=[x^ny^k]S(x,y)$, which denotes the coefficient of $x^ny^k$ in $S(x,y)$.
The first few values of this arrays are 
$$\mathcal{S}=[s(n,k)]_{n\geq 1, k\geq 0}=
\begin{pmatrix}
 1 & 0 & 0 & 0 & 0 \\
 2 & 0 & 0 & 0 & 0 \\
 5 & 0 & 0 & 0 & 0 \\
 13 & 1 & 0 & 0 & 0 \\
 35 & \framebox{\textbf{7}} & 0 & 0 & 0 \\
 98 & 33 & 1 & 0 & 0 \\
 284 & 135 & 10 & 0 & 0 \\
 846 & 519 & 64 & 1 & 0 \\
 2576 & 1933 & 340 & 13 & 0 
\end{pmatrix}.$$
For example, $s(5,1)=7$, the entry boxed in $\mathcal{S}$ above, and the corresponding Catalan words of length 5 with 1 symmetric valley are $$\tx{00101}, \quad \tx{01001}, \quad \tx{01010}, \quad  \tx{01011}, \quad  \tx{01012}, \quad \tx{01101}, \quad  \tx{01212}.$$ 
Let $s(n)$ be the sum of all symmetric valleys in the set of Catalan words of length $n$. 
\begin{coro}\label{corosym}
The generating function of the sequence
$s(n)$  is 
$$\sum_{n\geq 0}s(n)x^n=\frac{1 - 4 x + 2 x^2 - (1 - 2 x)\sqrt{1-4x}}{2(1-x)\sqrt{1-4x}}.$$
Moreover, for all $n\geq 1$, we have
\begin{align}
    s(n)=(3 n - 2)c_{n - 1}-\frac{1}{2}\sum_{k=1}^n\binom{2k}{k},\label{eq:cor sn}
\end{align}
 and an asymptotic approximation for $s(n)$ is $$\frac{2^{2n}}{12\sqrt{\pi n}}.$$
\end{coro}

For $n\geq 4$, the first few values of the sequence $s(n)$ are 
$$1, \quad 7, \quad 35, \quad 155, \quad 650, \quad 2652, \quad 10660, \quad 42484, \quad 168454, \quad  665874,\dots .$$
This sequence does not appear in the OEIS. 
An asymptotic approximation for the ratio between $s(n)$ and $v(n)$ is $1/2$.

\begin{openproblem}
    Find a combinatorial proof for the  expression for $s(n)$ given in equation \eqref{eq:cor sn}, as stated in Corollary~\ref{corosym}. 
\end{openproblem}

\section{The Distribution of Peaks}

\subsection{Peaks}
It is clear that the number of weak valleys of a Catalan word $w$ is exactly the number of weak peaks of $w$ minus one, but there is no natural link between the number of $\ell$-valleys and the number of $\ell$-peaks.
So, in order to count Catalan words according to the length and the number of $\ell$-peaks, we introduce the following bivariate generating function $$P_\ell(x,y)=\sum_{w \in\Cat^+}x^{|w|}y^{\lpea(w)}=\sum_{n\geq 1}x^{|w|}\sum_{w\in\Cat_n}y^{\lpea(w)},$$
where $\lpea(w)$ denotes the number of occurrences of subwords of the form $a(a+1)^\ell b$, and $a\geq b$, in $w$.  The coefficient of $x^ny^k$ in $P_\ell (x,y)$ is the number of Catalan words of length $n$ with $k$ $\ell$-peaks.

In Theorem~\ref{teosp1l}, we give an expression for this generating function.
\begin{theorem}\label{teosp1l}
The generating function for nonempty Catalan words with respect to the length and the  number of $\ell$-peaks is
$$P_\ell(x,y)={\frac {1-2x+x^{\ell+1}(1-y)-
\sqrt{1 - 4 x + 2 x^{\ell+ 1}(1 - y) + x^{2 (\ell+1)} (1 - y)^2}}{2({x}^{\ell+1}(y-1)+
x)}}.
$$
\end{theorem}
\begin{proof} From a similar argument as in the proof of  Theorem~\ref{teosv1l}, we have the functional equation:
\begin{align*}
P_\ell(x,y)=x \,+\, &  2xP_\ell(x,y)+\frac{x}{1-x}(P_\ell(x,y)-2xP_\ell(x,y)-x)P_\ell(x,y) \\
&+ xP_\ell(x,y)\left(\frac{x}{1-x}-x^\ell\right)+ xyx^\ell P_\ell(x,y)\\
&+x\left(P_\ell(x,y)\frac{x}{1-x}-P_\ell(x,y)x^\ell\right)P_\ell(x,y)+xyP_\ell(x,y)^2x^\ell.
\end{align*}
Solving this equation yields the desired result.
\end{proof}
Let $p_\ell(n)$ be the sum of all $\ell$-peaks in the set of Catalan words of length $n$. 
\begin{coro}
The generating function of the sequence
$p_\ell(n)$  is 
$$\sum_{n\geq 1}p_\ell(n)x^n={\frac {{x}^{\ell-1} \left( 1-3\,x-(1-x)\sqrt {1-4\,x}
 \right) }{2\sqrt {1-4\,x}}}.$$
Moreover, for all $n\geq 1$ we have
\begin{align}\label{eq:pell}
p_\ell(n)&=\binom{2(n - \ell) - 1}{n - \ell - 2},
\end{align}
 and an asymptotic approximation for $p_\ell(n)$ is $$\frac{2^{2(n-\ell)-1}}{\sqrt{\pi n}}.$$
\end{coro}

Taking $\ell=1$ in Theorem~\ref{teosp1l}, establishes that the generating function for Catalan words with respect to the length and the  number of short peaks is
$$P_1(x,y)={\frac {1-2x+x^2(1-y)-
\sqrt{1 - 4 x + 2 x^2(1 - y) + x^4 (1 - y)^2}}{2({x}^2(y-1)+
x)}}.
$$

Let $p_1(n,k)$ denote the number of Catalan words of length $n$ with exactly $k$ short peaks, that is $p_1(n,k)=[x^ny^k]P_1(x,y)$, which denotes the coefficient of $x^ny^k$ in $P_1(x,y)$. The first few values of this array are

$$\mathcal{P}_1=[p_1(n,k)]_{n\geq 1, k\geq 0}=
\begin{pmatrix}
 1 & 0 & 0 & 0 &0\\
 2 & 0 & 0 & 0 &0\\
 4 & 1 & 0 & 0& 0\\
 9 & 5 & 0 & 0 &0\\
 22 & 19 & 1 & 0&0 \\
 57 & 66 & 9 & 0 &0\\
 154 & 221 & 53 & 1 &0\\
 429 & 729 & 258 & 14 &0\\
 1223 & 2391 & 1131 & 116&1 \\
\end{pmatrix}.$$
The array $\mathcal{P}_1$ corresponds to the OEIS entry \cite[\seqnum{A116424}]{OEIS}. Notice that in the notation of subwords patterns, see  \cite{Shattuck}, the sequence $p_1(n)$ counts the  joint distribution of the  patterns $121$ and $231$.

Using a similar proof as for Theorem~\ref{teosp1l}, we generalize the result in order to obtain the following generating function for the
number Catalan words of length $n$ with respect to the number of peaks.

\begin{theorem}\label{teop1}
The generating function for Catalan words with respect to the length and the  number of peaks is
$$P(x,y)=\scalemath{0.95}{{\frac {1-3x-{x}^{2}y+3\,{x}^{2}-\sqrt {1 - 6 x + 11 x^2 - 6 x^3 + x^4 - 2 x^2 y + 2 x^3 y - 2 x^4 y + x^4 y^2}}{2x
 \left( xy-2\,x+1 \right) }}}
.$$
\end{theorem}

\begin{figure}[H]
\centering
\includegraphics[scale=0.8]{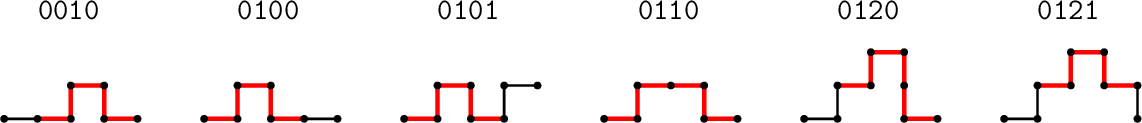}
\caption{Catalan words of length 4 with one peak.} \label{ExPeak}
\end{figure}

Let $p(n,k)$ denote the number of Catalan words of length $n$ with exactly $k$ peaks, that is $p(n,k)=[x^ny^k]P(x,y)$, which denotes the coefficient of $x^ny^k$ in $P(x,y)$.
The first few values of this arrays are 
$$\mathcal{P}=[p(n,k)]_{n\geq 1, k\geq 0}=
\begin{pmatrix} 
 1 & 0 & 0 & 0 & 0 \\
 2 & 0 & 0 & 0 & 0 \\
 4 & 1 & 0 & 0 & 0 \\
 8 &  \framebox{\textbf{6}} & 0 & 0 & 0 \\
 16 & 25 & 1 & 0 & 0 \\
 32 & 89 & 11 & 0 & 0 \\
 64 & 290 & 74 & 1 & 0 \\
 128 & 893 & 392 & 17 & 0 \\
 256 & 2645 & 1796 & 164 & 1 \\
\end{pmatrix}.$$
This array does not appear in the OEIS. 
For example, $p(4,1)=6$, the entry boxed in $\mathcal{P}$ above, and the corresponding Catalan words of length 4 with one peak  are shown in Figure~\ref{ExPeak}.

Let $p(n)$ be the sum of all peaks in the set of Catalan words of length $n$. 
\begin{coro}
The generating function of the sequence
$p(n)$  is 
$$\sum_{n\geq 0}p(n)x^n=\frac {1-3\,x-(1-x)\sqrt {1- 4\,x}}{2(1-x)\sqrt {1-4x}}
,$$
 and an asymptotic approximation for $p(n)$ is  $$\frac{2^{2n}}{6\sqrt{\pi n}}.$$
\end{coro}

The first few values of the sequence $p(n)$ ($n\geq 3$) are 
$$1, \quad 6, \quad 27, \quad 111, \quad 441, \quad 1728, \quad  6733, \quad  26181,\quad  101763, \quad 395693, \dots .$$
This sequence does not appear in the OEIS. 
Since $p(n)=\displaystyle\sum_{\ell\geq 1}p_\ell(n)$, by equation \eqref{eq:pell}, we deduce 
$$p(n)=\sum_{\ell=1}^{n-1}\binom{2(n - \ell) - 1}{n - \ell - 2}.$$

\subsection{Symmetric Peaks}
A \emph{symmetric peak} is a peak of the form $a(a+1)^\ell a$ with $\ell\geq 1$. 
Let $\psym(w)$ denote  the number of the symmetric peaks of the word $w$.  
In order to count Catalan words according to the length and the number symmetric peaks, we introduce the following bivariate generating function
$$T(x,y)=\sum_{w\in\Cat^+}x^{|w|}y^{\psym(w)}=\sum_{n\geq 1}x^{|w|}\sum_{w\in\Cat_n}y^{\psym(w)},$$
where the coefficient of $x^ny^k$ in $T(x, y)$ is the number of Catalan words of length $n$ with $k$ symmetric peaks.

In Theorem~\ref{thm:gf sympeaks}, we give an expression for this generating function.
\begin{theorem}\label{thm:gf sympeaks}
    The generating function of the nonempty Catalan words with respect to the length and the number of symmetric peaks is 
    $$T(x,y)=\scalemath{0.92}{{\frac {{x}^{2}y-3\,{x}^{2}+3\,x-1+\sqrt {{x}^{4}{y}^{2}-6\,{x}^{
4}y+5\,{x}^{4}+6\,{x}^{3}y-10\,{x}^{3}-2\,{x}^{2}y+11\,{x}^{2}-6\,x+1}
}{2x \left( -1+x \right) }}}
.$$    
\end{theorem}
\begin{proof} Let $w$ be a nonempty Catalan word, and let  $w=\texttt{0}(w'+1)w''$ be the first return decomposition, with $w', w''\in \Cat$. If $w'=w''=\epsilon$, then $w=\texttt{0}$, and its generating function is $x$. If $w'\neq \epsilon$ and $w''=\epsilon$, then $w=\texttt{0}(w'+1)$, and its generating function is $xT(x,y)$. Similarly,  if $w'=\epsilon$ and $w''\neq \epsilon$, then $w=\texttt{0}w''$, and its generating function is $xT(x,y)$.  Finally, if $w'$ and $w''$ are not empty, then $w=\texttt{0}(w'+1)w''$, and we have a new symmetric peak in the case where $w'$ is $\tx{0}^k$, $k\geq 1$.  The generating for this case is
\[\frac{x^2y}{1-x}T(x,y)+x\left(T(x,y)-\frac{x}{1-x}\right)T(x,y).\]
Summing  all the cases and solving the obtained functional equation yields the desired result.
\end{proof}

Let $t(n,k)$ denote the number of Catalan words of length $n$ with exactly $k$ symmetric peaks, that is $t(n,k)=[x^ny^k]T(x,y)$, which denotes the coefficient of $x^ny^k$ in $T(x,y)$.
The first few values of this arrays are 
$$\mathcal{T}=[t(n,k)]_{n\geq 1, k\geq 0}=
\begin{pmatrix}
 1 & 0 & 0 & 0 & 0 \\
 2 & 0 & 0 & 0 & 0 \\
 4 & 1 & 0 & 0 & 0 \\
 9 & \framebox{\textbf{5}}  & 0 & 0 & 0 \\
 23 & 18 & 1 & 0 & 0 \\
 64 & 60 & 8 & 0 & 0 \\
 187 & 199 & 42 & 1 & 0 \\
 563 & 667 & 189 & 11 & 0 \\
 1731 & 2259 & 795 & 76 & 1 
\end{pmatrix}.$$

For example, $t(4,1)=5$, the entry boxed in $\mathcal{T}$ above, and the corresponding Catalan words of length 4 with 1 symmetric peak are
$$\tx{0100}, \quad \tx{0101}, \quad \tx{0010}, \quad  \tx{0110},\quad  \tx{0121}.$$

Let $t(n)$ be the sum of all symmetric peaks in the set of Catalan words of length $n$. 
\begin{coro}
The generating function of the sequence
$t(n)$  is 
$$\sum_{n\geq 0}t(n)x^n={\frac {x \left( 2\,{x}^{2}-3\,x+1-(1-x)\sqrt {1-4x} \right) }{2\left( 1-x \right)^2\sqrt {1-4x}  }}
.$$
Moreover, for all $n\geq 1$,  we have
$$t(n)=\sum_{k=0}^{n-3}\binom{2k+2}{k},$$
 and an asymptotic approximation for $t(n)$ is  $$\frac{2^{2n}}{12\sqrt{\pi n}}.$$
\end{coro}

For $n\geq 3$, the first few values of the sequence $t(n)$ are 
$$1, \quad 5, \quad 20, \quad 76, \quad 286, \quad 1078, \quad 4081, \quad 15521, \quad 59279, \quad  227239, \dots .$$
This sequence corresponds to the OEIS entry  \cite[\seqnum{A057552}]{OEIS}. 
An asymptotic approximation for the ratio between $t(n)$ and $p(n)$ is $1/2$.\\

\noindent{\bf Acknowledgement:}
Pamela E.~Harris was supported in part by a Karen Uhlenbeck EDGE Fellowship.
Jos\'e L.~Ramírez was partially supported by Universidad Nacional de Colombia. The last author  would like to thank the Laboratoire d’Informatique de Bourgogne for the warm hospitality during his visit at Universit\'e de Bourgogne where part of this work was done.

\end{document}